# Orlicz regrets to consistently bound statistics of random variables with an application to environmental indicators

Hidekazu Yoshioka[a]* and Yumi Yoshioka[b]

[a]*Japan Advanced Institute of Science and Technology, 1-1 Asahidai, Nomi, Ishikawa 923-1292, Japan*

[b]*Gifu University, 1-1 Yanagido, Gifu, Gifu 501-1193, Japan*

* Corresponding author: yoshih@jaist.ac.jp, ORCID: 0000-0002-5293-3246

**Abstract** Evaluating environmental variables that vary stochastically is the principal topic for designing better environmental management and restoration schemes. Both the upper and lower estimates of these variables, such as water quality indices and flood and drought water levels, are important and should be consistently evaluated within a unified mathematical framework. We propose a novel pair of Orlicz regrets to consistently bound the statistics of random variables both from below and above. Here, consistency indicates that the upper and lower bounds are evaluated with common coefficients and parameter values being different from some of the risk measures proposed thus far. Orlicz regrets can flexibly evaluate the statistics of random variables based on their tail behaviour. The explicit linkage between Orlicz regrets and divergence risk measures was exploited to better comprehend them. We obtain sufficient conditions to pose the Orlicz regrets as well as divergence risk measures, and further provide gradient descent-type numerical algorithms to compute them. Finally, we apply the proposed mathematical framework to the statistical evaluation of water quality data as key environmental indicators in a Japanese river environment.

**Competing interests** The authors declare no competing interests.

**Funding** Japan Society for the Promotion of Science No. 22K14441 and 22H02456.

**Acknowledgements** The authors thank Dr. Ikuo Takeda of Shimane University for kindly providing time-series data on the water quality indices at Otsu in the midstream Hii River.

**Availability of data and material** Data will be available upon reasonable request from the corresponding author.

**Author contributions** *Hidekazu Yoshioka*: Conceptualization, Methodology, Software, Formal analysis, Data Curation, Visualization, Writing original draft preparation, writing review and editing, Supervision, Project administration, funding acquisition; *Yumi Yoshioka*: Data Curation, Visualization, Writing original draft preparation, writing review, and editing.

**Declaration** No AI technologies were used throughout this study.

# 1. Introduction

## 1.1 Research background

The sustainable coexistence between human society and the environment is an urgent global issue (Gao and Clark, 2023; Henderson and Loreau, 2023; Yin et al., 2023). Designing effective schemes for environmental management and restoration is pivotal for addressing this issue. Environmental variables, such as the quantity and quality of surface water (Li et al., 2023; Piffer et al., 2023), amount of rainfall (Ulzega and Albert, 2023), and quality of air (Perri rt al., 2023), change stochastically. A mathematical framework that can address these characteristics should be used to assess the environmental variables.

Stochastic models are essential mathematical tools for formulating and analysing environmental variables. Perri and Porporato (2022) studied a jump-driven model of soil moisture and salinity and fitted the Beta distribution to salinity concentration. Tu et al. (2023) mapped the noise colour, that is, the spectral shape, of discharges over river networks on a continental scale. Park and Kim (2023) proposed a discrete-time stochastic model to simulate stochastic snow depth time series in mountainous regions. Sumata et al. (2023) studied the regime shifts in ocean sea ice driven by stochastic climate change. Larsson (2023) proposed parametric heatwave insurance, a type of derivative based on the option pricing theory of stochastic processes.

Conventional stochastic models have a common drawback in that their parameters and coefficients are often conceptual, and hence, lumped ones that will be subject to measurement as well as structural errors (Mészáros et al., 2021; Oberpriller et al., 2021). They can be estimated from real data using statistical (Smith et al., 2023; Yoshioka and Yamazaki, 2023), data-driven (Nikakhtar et al., 2023; Yang et al., 2023), or Bayesian methods (Reichert et al., 2021), while they are not free from estimation errors. Therefore, studies based on a stochastic model should include some uncertainty analysis so that model performance under the errors can be discussed.

Risk measures to reasonably overestimate extreme events that will be harmful for decision-makers have been developed in the finance and economics sector (Ben-Tal and Teboulle, 2007; Chan and Hu, 2018). A major class of risk measures are based on statistical divergences (Pichler, and Schlotter, 2020): convex functions to evaluate the proximity between two stochastic models (Sason and Verdú, 2016). The widely-used divergences are Kullback–Leibler divergence also called relative entropy and its generalization such as $\alpha$, Rényi, Tsallis, and functional divergences (Van Erven and Harremos, 2014; Birrell et al., 2023; Nikoufar and Kanani Arpatapeh, 2023). Sophisticated risk measures include the entropic value-at-risk and its extended versions (Ahmadi-Javid, 2012; Bodnar et al., 2022; Zou et al., 2023), and several coherent risk measures (Eshghali and Krokhmal, 2023) have also been proposed. These risk measures have typically been applied to the overestimation problem of a target random variable, while the underestimation problem has been addressed less frequently, except in a few recent studies (Christensen and Connault, 2023; Yoshioka and Yoshioka, 2023a; Yoshioka and Yoshioka, 2023b; Yoshioka et al., 2023b).

Both overestimation and underestimation problems are important for environmental research applications. For example, water levels that are too high or too low, which represent floods and droughts, respectively, are of high importance. Both the upper and lower bounds of the target statistic must be evaluated efficiently and consistently. Here, efficiency implies that the quantity can be computed using an analytical calculation or a tractable numerical method. Consistency implies that the upper and lower bounds are defined within a unified setting (e.g., using the same coefficients and parameters). Employing different methodologies for the upper and lower bounds potentially leads to biased results. The authors recently proposed statistical methodologies to bound environmental variables both from below and above, but they are inconsistent in the above-mentioned sense because different shapes of risk measures were used to evaluate the two bounds (Yoshioka and Yoshioka, 2023a; Yoshioka and Yoshioka, 2023b; Yoshioka et al., 2023b); the authors did not have the idea for finding a consistent framework. This limitation must be overcome to establish an efficient and theoretically better-consistent method for assessing the environment.

### 1.2 Objective and contribution

The goal of this study is to formulate and analyse a novel pair of risk measures to bound a random variable both from below and above, so that the upper and lower bounds of its statistics are evaluated consistently. The core of our method is the recently proposed Orlicz regret measure (Fröhlich and Williamson, 2023), called Orlicz regret, which was originally proposed to evaluate the upper bound of the expectation of random variables in the context of machine learning. The name "Orlicz" comes from Orlicz spaces as a generalization of Lebesgue spaces (Rao and Ren, 2002), which are functional spaces to analyse random variables under general growth conditions. Various statistics, including polynomial and exponential moments, can be evaluated in detail as elements belonging to Orlicz spaces. It should be noted that risk measures can be systematically studied through suitable Orlicz spaces (Cheridito and Li, 2009; Bellini et al., 2018; Liu and Shushi, 2023). In this study, we extend the theory of Fröhlich and Williamson (2023) to the upper bound of random variables such that the lower bound is also consistently covered.

A remarkable property of the Orlicz regret is that it has an explicit linkage with constrained maximisation and minimisation problems subject to a generic statistical divergence: the divergence risk measure. More specifically, not only do Orlicz regrets serve as risk measures by themselves, but the associated divergence risk measures can also be expressed through suitable inf- and sup- convolutions of Orlicz regrets (Vinel and Krokhmal, 2017). This provides a clear statistical interpretation of the Orlicz regret. The convex conjugate of divergences is effectively used to better comprehend Young's functions and the Orlicz regret. Specifically, the convex conjugate of the divergence serves as Young's function in the corresponding Orlicz space to which the target random variable belongs. Another remarkable property of Orlicz regret is that it is based on the equivalence between Orlicz and Amemiya norms (Hudzik and Maligranda, 2000), through which we show that both the upper- and lower-bounding Orlicz regrets can be formulated consistently using common parameters and coefficients as desired. Moreover, as a by-product, we can construct a simple gradient descent-type method to numerically compute the Orlicz regret pair and the associated divergence risk measures.

Finally, we apply the proposed mathematical framework to unique water-quality data collected from rivers in Japan. The concentrations of water quality indices were fitted to a gamma distribution; for the gamma distribution, we obtained a sharp, well-posed condition for the Orlicz regrets. Their parameter sensitivities were numerically analysed along with added computational results concerning the safety probability against the eutrophication as an extreme event.

The rest of this paper is organized as follows. The first half of **Section 2** presents some preliminaries to be used in this study and reviews the Orlicz regret for the upper bound. This section proposes a mutually consistent lower-bound counterpart. The latter half of **Section 2** is devoted to analysing the Orlicz regrets and the associated divergence risk measures, focusing on their well-posedness and the gamma case. **Section 3** addresses the application of the Orlicz regrets to the statistical evaluation of environmental variables. **Section 4** concludes the paper and presents future perspectives. **Appendices** contain technical proofs of the propositions stated in the main text and some technical result.

## 2. Orlicz regrets

### 2.1 Preliminaries

Several notions used in this study are introduced in this subsection. The explanation below is based on the literature (Fröhlich and Williamson, 2023; Rao and Ren, 2002; Hudzik and Maligranda, 2000). We consider a real-valued scalar random variable denoted as $X$ in the complete probability space $(\Omega, \mathcal{F}, \mathbb{P})$. For a generic probability measure $\mathbb{Q}$, the expectation is denoted as $\mathbb{E}_{\mathbb{Q}}[\cdot]$. If there is no confusion, the expectation is simply expressed as $\mathbb{E}[\cdot]$. We always assume that $X$ is measurable with respect to $\mathbb{P}$ and that all (in)equalities of the random variables are almost surely satisfied.

The extension of any function $F:[0,+\infty)\to[0,+\infty)\cup\{+\infty\}$ to $\mathbb{R}$ is defined such that $F(x)=+\infty$ if $x<0$. The Young's function $\Phi:[0,+\infty)\to[0,+\infty)\cup\{+\infty\}$ is a non-negative, increasing, and convex function that satisfies $\Phi(0)=0$. For any convex function $F:[0,+\infty)\to[0,+\infty)\cup\{+\infty\}$, the convex conjugate $F_c:\mathbb{R}\to\mathbb{R}\cup\{+\infty\}$ is defined as $F_c(y)=\sup_{x\in\mathbb{R}}\{xy-F(x)\}$. The convex conjugate of the Young function $\Phi$ is denoted as $\Psi$, the latter being a Young function as well.

The statistical divergence, referred as divergence below for simplicity, is a proper, lower-semicontinuous, and convex function $f:[0,+\infty)\to[0,+\infty)\cup\{+\infty\}$ such that $f(1)=0$, $f(0)<+\infty$, $f(x)$ is bounded for any $x\geq 0$, and is twice continuously differentiable at any $x>0$. The convex conjugate of $f$ is denoted as $g$. With $\varepsilon>0$, set $f_\varepsilon=\varepsilon^{-1}f$. The parameter $\varepsilon>0$ serves as the uncertainty aversion level as shown in the later subsections. The convex conjugate of $f_\varepsilon=\varepsilon^{-1}f$ is given by $g_\varepsilon=\varepsilon^{-1}g(\varepsilon\cdot)$. Finally, set functions related to Orlicz regrets:

$$\bar{g}(y)=\begin{cases}0 & (y<0)\\ g(y) & (y\geq 0)\end{cases}. \tag{1}$$

This $\bar{g}$ serves as a Young's function, and is the convex conjugate of another Young's function

$$\bar{f}(x)=\begin{cases}0 & (x<1)\\ f(x) & (x\geq 1)\end{cases}. \tag{2}$$

***Remark 1*** The requirement $f(x)=+\infty$ for $x<0$ appears unusual, but it is essential to obtain the representation formula of the worst-case estimation of **Proposition 2** presented later.

***Remark 2*** The function $f$ is expressed as $f(x) = x\ln x - x + 1$ for the Kullback–Leibler divergence ($0\ln 0 = 0$) with the conjugate $g(y) = e^y - 1$. For the $\alpha$ divergence with $\alpha > 0$, we have $f(x) = \frac{x^\alpha - \alpha x + \alpha - 1}{\alpha - 1}$. An elementary calculation yields $g(y) = \left(1 + \frac{\alpha - 1}{\alpha} y\right)_+^{\frac{\alpha}{\alpha-1}} - 1$ ($(\cdot)_+ = \max\{0,\cdot\}$), which is bounded only for $1 + \frac{\alpha-1}{\alpha} y > 0$ when $\alpha \in (0,1)$ (**Figure 1**). The $\alpha$ divergence is equivalent to the Kullback–Leibler one under the limit $\alpha \to 1$.

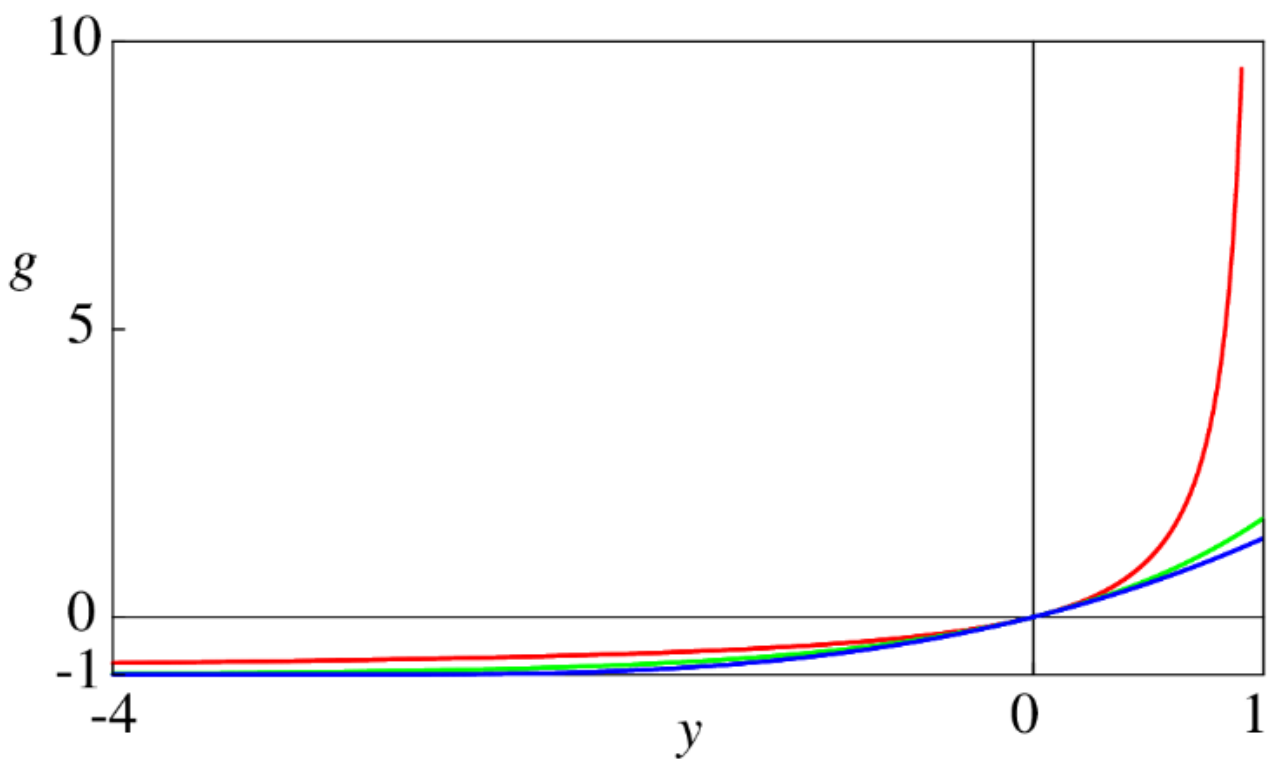

**Figure 1.** Graphs of $g$ with $\alpha = 0.5$ (red, increases exponentially at $y = 1$), $\alpha = 1.0$ (green), and $\alpha = 1.5$ (blue, equals $-1$ for $y \le -3$).

We define a generic Orlicz space and its associated norms. We fix Young's function $\Phi$ with the convex conjugate $\Psi$. An Orlicz space is a Banach space as a collection of random variables $X$ such that the following functional space, called Orlicz space, is nonempty:

$$\mathcal{L}_\Phi = \left\{ X \middle| \mathbb{E}\left[\Phi\left(\frac{|X|}{\lambda}\right)\right] < +\infty \text{ for some } \lambda > 0 \right\}. \tag{3}$$

The Orlicz space is a Banach space equipped with a norm called the Luxemburg norm:

$$\|X\|_{\Phi,L} = \inf\left\{ \lambda > 0 \middle| \mathbb{E}\left[\Phi\left(\frac{|X|}{\lambda}\right)\right] \le 1 \right\}. \tag{4}$$

There exist two norms equivalent to the Luxemburg norm, which are the Orlicz norm

$$\|X\|_{\Phi,O} = \sup\left\{ \mathbb{E}[XZ] \mid \mathbb{E}[|Z|] < +\infty,\ \mathbb{E}[\Psi(|Z|)] \le 1 \right\} \tag{5}$$

and Amemiya norm, the latter being linked to the Orlicz regrets:

$$\|X\|_{\Phi,A} = \inf_{t>0}\left\{ t\left(1 + \mathbb{E}\left[\Phi\left(\frac{|X|}{t}\right)\right]\right) \right\}. \tag{6}$$

In (5) defining the Orlicz norm, the measurable random variable $Z$ plays a role of a Radon–Nikodým derivative between two probability measures. According to Theorem 1 of Hudzik and Maligranda (2000), the Orlicz and Amemiya norms are equal; for any $X \in \mathcal{L}_{\Phi}$, it follows that $\|X\|_{\Phi,O} = \|X\|_{\Phi,A}$.

***Remark 3*** It is shown in the subsequent subsections that $\Psi$ is formally observed as a divergence $f$ while $\Phi$ as the convex conjugate $g$. This perspective is useful in both theory and application as it implies that the Young function and divergence, which have been specified separately in the previous studies (Yoshioka et al., 2023b; Bellini et al., 2018), can be specified in a unified manner. This is a theoretical advantage of the Orlicz regrets as demonstrated in this paper. Its potential disadvantage is the lack of the translation invariance, whereas divergence risk measures ((8) and (11)) satisfy it. This is actually not a significant drawback in our application, while it may become crucial for some applications such as assessing monetary risks.

## 2.2 Upper-bound

In this study, we extend the theory of Fröhlich and Williamson (2023) to both upper- and lower-bounding Orlicz regrets. We fix a divergence $f$. For any $X \in \mathcal{L}_{\bar{g}}$, the upper-bounding Orlicz regret $V_g$ is given by

$$V_g(X) = \inf_{t>0}\left\{t\left(1+\mathbb{E}\left[g\left(\frac{X}{t}\right)\right]\right)\right\} = \sup\left\{\mathbb{E}[XZ] \mid Z \in \mathcal{L}_{\bar{f}},\ Z \geq 0,\ \mathbb{E}\left[f(Z)\right] \leq 1\right\}. \tag{7}$$

We have the equivalence between the Amemiya norm $V_g(X) = \|X\|_{g,A}$ if $X \geq 0$, which is why $V_g$ is called the "Orlicz" regret. The restriction $Z \in \mathcal{L}_{\bar{f}}$ leads to a feasible solution on the right-hand side of (7) because if it is not true, $\mathbb{E}\left[f(Z)\right]$ must be $\mathbb{E}\left[f(Z)\right] = +\infty$ due to the profile of $f$, which is a contradiction.

Orlicz regret serves as an asymmetric norm (and hence, it is not a conventional norm), measuring both the gain ($X > 0$) and loss ($X < 0$) in a unified manner; here, the former corresponds to the overestimation while the latter to the underestimation, respectively. This $V_g$ has the following properties (Proposition 3.9 of Fröhlich and Williamson (2023)):

**V1.** Positive homogeneity: $V_g(\lambda X) = \lambda V_g(X)$ for any $\lambda \geq 0$.

**V2.** Subadditivity: $V_g(X+Y) \leq V_g(X) + V_g(Y)$.

**V3.** Monotonicity: If $X \leq Y$, then $V_g(X) \leq V_g(Y)$.

**V4.** Aversity: $\mathbb{E}[X] \leq V_g(X)$.

**V5.** Point-separating: $V_g(X) = 0$ if and only if $X = 0$.

We now define the upper-bounding divergence risk measure, which serves as the constrained worst-case overestimation of the expectation of $X$ under uncertainty:

$$\bar{R}_{f,\varepsilon}(X)=\sup\left\{\mathbb{E}_{\mathbb{Q}}[X] \mid \mathbb{Q}\in\mathfrak{Q},\ \mathbb{E}_{\mathbb{P}}\left[f\left(\frac{\mathrm{d}\mathbb{Q}}{\mathrm{d}\mathbb{P}}\right)\right]\le\varepsilon\right\}=\sup\left\{\mathbb{E}_{\mathbb{Q}}[X] \mid \mathbb{Q}\in\mathfrak{Q},\ \mathbb{E}_{\mathbb{P}}\left[f_{\varepsilon}\left(\frac{\mathrm{d}\mathbb{Q}}{\mathrm{d}\mathbb{P}}\right)\right]\le 1\right\}, \tag{8}$$

where $\frac{\mathrm{d}\mathbb{Q}}{\mathrm{d}\mathbb{P}}$ denotes the Radon–Nikodým derivative of $\mathbb{Q}$ with respect to $\mathbb{P}$ and $\mathfrak{Q}$ denotes the set of probability measures absolutely continuous with respect to $\mathbb{P}$. According to Fröhlich and Williamson (2023), $\bar{R}_{f,\varepsilon}(X)$ is expressed as the inf-convolution of $V_{g_\varepsilon}$:

$$\bar{R}_{f,\varepsilon}(X)=\inf_{\mu\in\mathbb{R}}\left\{\mu+V_{g_\varepsilon}(X-\mu)\right\}=\inf_{\substack{\mu\in\mathbb{R}\\ t>0}}\left\{\mu+t\left(1+\mathbb{E}\left[g_{\varepsilon}\left(\frac{X-\mu}{t}\right)\right]\right)\right\}. \tag{9}$$

In the later section, we show that an analogous formulation holds true for the lower bounding counterpart. This representation formula explicitly connects the divergence risk measure and Young's function through the Orlicz regret.

### 2.3 Lower-bound

The lower-bounding Orlicz regret is proposed in this study, which is dual to the upper-bounding regret proposed in the earlier subsection. In addition, we obtain a dual representation of the lower-bounding divergence risk measure. We fix a divergence $f$. For any $X\in\mathcal{L}_{\bar{g}}$, following the upper-bounding case, we define the lower-bounding Orlicz regret as

$$W_g(X)=\sup_{t>0}\left\{-t\left(1+\mathbb{E}\left[g\left(\frac{-X}{t}\right)\right]\right)\right\}=\inf\left\{\mathbb{E}[XZ] \mid Z\in\mathcal{L}_{\bar{f}},\ Z\ge 0,\ \mathbb{E}\left[f(Z)\right]\le 1\right\}. \tag{10}$$

Its properties are summarized in **Proposition 1** below, which are dual to **V1-V5**.

***Proposition 1***

*The Orlicz regret $W_g$ satisfies the following properties **W1-W5**. Furthermore, the right equality of (10) holds true.*

***W1.*** *Positive homogeneity:* $W_g(\lambda X)=\lambda W_g(X)$ *for any* $\lambda\ge 0$.

***W2.*** *Superadditivity:* $W_g(X+Y)\ge W_g(X)+W_g(Y)$.

***W3.*** *Monotonicity: If* $X\le Y$, *then* $W_g(X)\le W_g(Y)$.

***W4.*** *Aversity:* $\mathbb{E}[X]\ge W_g(X)$.

***W5.*** *Point-separating: If* $X\le 0$, *then* $W_g(X)=0$ *if and only if* $X=0$.

We now define the lower-bounding divergence risk measure, which serves as the constrained worst-case underestimation of the expectation of $X$ under uncertainty:

$$\underline{R}_{f,\varepsilon}(X)=\inf\left\{\mathbb{E}_{\mathbb{Q}}[X] \mid \mathbb{Q}\in\mathfrak{Q},\ \mathbb{E}_{\mathbb{P}}\left[f\left(\frac{\mathrm{d}\mathbb{Q}}{\mathrm{d}\mathbb{P}}\right)\right]\leq\varepsilon\right\}=\inf\left\{\mathbb{E}_{\mathbb{Q}}[X] \mid \mathbb{Q}\in\mathfrak{Q},\ \mathbb{E}_{\mathbb{P}}\left[f_{\varepsilon}\left(\frac{\mathrm{d}\mathbb{Q}}{\mathrm{d}\mathbb{P}}\right)\right]\leq 1\right\}. \quad (11)$$

Its representation through Orlicz regret has not been found in the previous studies, although it is reasonable to expect that some Orlicz regret can be used to define $\underline{R}_{f,\varepsilon}$. **Proposition 2** shows that $W_{g_\varepsilon}$ is the desired one. Consequently, we can formulate a lower-bounding Orlicz regret consistent with the upper-bounding regret because they share the same coefficients.

***Proposition 2***

*It follows that*

$$\underline{R}_{f,\varepsilon}(X)=\sup_{\mu\in\mathbb{R}}\left\{\mu+W_{g_\varepsilon}(X-\mu)\right\}=\sup_{\substack{\mu\in\mathbb{R}\\ t>0}}\left\{\mu-t\left(1+\mathbb{E}_{\mathbb{P}}\left[g_{\varepsilon}\left(-\frac{X-\mu}{t}\right)\right]\right)\right\}. \quad (12)$$

*That is,* $\underline{R}_{f,\varepsilon}$ *is a sup-convolution of* $W_{g_\varepsilon}$.

## 2.4 Well-posedness results

For obtaining the meaningful Orlicz regret and associated divergence risk measure such that they do not diverge to $\pm\infty$, we must ensure that the expectation $\mathbb{E}\left[g_{\varepsilon}\left(\pm\frac{X-\mu}{t}\right)\right]$ exists as a real value, which may not be the case for the $\alpha$ divergence with $\alpha\in(0,1)$ because of

$$g_{\varepsilon}(y)=\frac{1}{\varepsilon}\left(1-\frac{1-\alpha}{\alpha}\varepsilon y\right)^{\frac{-\alpha}{1-\alpha}}-1,\ y<\frac{\alpha}{(1-\alpha)\varepsilon} \text{ and } g_{\varepsilon}(y)=+\infty \text{ otherwise.} \quad (13)$$

Thereby, this $g_\varepsilon$ is not globally defined for $y\in\mathbb{R}$ and diverges to $+\infty$ at the positive $y$.

We analyse this issue in detail for non-negative and unbounded $X$ (e.g., gamma random variables) in **Proposition 3** below; the proposition suggests that this singularity does not deteriorate the well-posedness of the Orlicz regret and worst-case estimate for the lower bound, whereas the Kullback–Leibler case $\alpha=1$ may not be robust against perturbations as for the upper bound because of the non-existence of any minimising $(\mu,t)\in\mathbb{R}\times(0,+\infty)$ of $\bar{R}_{f,\varepsilon}$ for any $\alpha\in(0,1)$. The $\alpha$ divergence with $\alpha>1$ is more robust than the Kullback–Leibler case in this sense. **Proposition 3** shows that $\bar{R}_{f,\varepsilon}$ and $\underline{R}_{f,\varepsilon}$ can be defined for the $\alpha$ divergence if the range of $(\mu,t)$ is suitably restricted for appropriate $\alpha$, which will be assumed below. Specifically, $\alpha>1$ for the upper bound, unless otherwise specified.

***Proposition 3***

*Assume that* $X\geq 0$ *and* $X$ *has a probability density function that is positive on* $(0,+\infty)$. *Assume further that* $\alpha\in(0,1)$. *Then, there is no minimising pair* $(\mu,t)\in\mathbb{R}\times(0,+\infty)$ *of* $\bar{R}_{f,\varepsilon}$. *Furthermore, any maximising pair* $(\mu,t)\in\mathbb{R}\times(0,+\infty)$ *of* $\underline{R}_{f,\varepsilon}$ *satisfies* $\frac{\alpha}{(1-\alpha)\varepsilon}\geq\frac{\mu}{t}$.

### 2.5 Ordering properties

There is an important ordering property between Orlicz regrets and divergence risk measures. **Proposition 4** shows that both the pairs $\left(\underline{R}_{g,\varepsilon},\overline{R}_{g,\varepsilon}\right)$ and $\left(W_{g_\varepsilon},V_{g_\varepsilon}\right)$ can bound the target expectation $\mathbb{E}_{\mathbb{P}}\left[X\right]$ for each uncertainty aversion level $\varepsilon>0$, and the latter pair is more pessimistic. The degree of pessimism is analysed numerically later.

***Proposition 4***

*For any* $\varepsilon>0$, *it follows that*

$$W_{g_\varepsilon}\left(X\right)\le \underline{R}_{f,\varepsilon}\left(X\right)\le \mathbb{E}_{\mathbb{P}}\left[X\right]\le \overline{R}_{f,\varepsilon}\left(X\right)\le V_{g_\varepsilon}\left(X\right). \tag{14}$$

**Proposition 5** shows the parameter dependence of Orlicz regret and the associated risk measures. These inequalities show the monotone dependence of regrets and risks against the uncertainty aversion parameter $\varepsilon$, such that a larger uncertainty aversion results in a more pessimistic estimate.

***Proposition 5***

*For any* $0<\varepsilon_1\le\varepsilon_2$, *it follows that*

$$\overline{R}_{f,\varepsilon_1}\left(X\right)\le \overline{R}_{f,\varepsilon_2}\left(X\right) \text{ and } \underline{R}_{f,\varepsilon_2}\left(X\right)\le \underline{R}_{f,\varepsilon_1}\left(X\right) \tag{15}$$

*as well as,*

$$V_{g_{\varepsilon_1}}\left(X\right)\le V_{g_{\varepsilon_2}}\left(X\right) \text{ and } W_{g_{\varepsilon_2}}\left(X\right)\le W_{g_{\varepsilon_1}}\left(X\right). \tag{16}$$

*It also follows that*

$$\frac{\varepsilon_2}{\varepsilon_1}V_{g_{\varepsilon_1}}\left(X\right)\le V_{g_{\varepsilon_2}}\left(X\right) \text{ for } X\ge 0 \text{ and } W_{g_{\varepsilon_2}}\left(X\right)\ge \frac{\varepsilon_2}{\varepsilon_1}W_{g_{\varepsilon_1}}\left(X\right) \text{ for } X\le 0. \tag{17}$$

***Remark 4*** The Radon–Nikodým derivatives $Z$ realizing the divergence risk measures can be obtained as follows. For the lower bound, assume that we have a minimising pair $\left(\mu,t\right)\in\mathbb{R}\times\left(0,+\infty\right)$. Subsequently, by solving the inner maximisation problem of (41) with respect to $Z$, we obtain the minimiser

$$Z=h_\varepsilon\left(-\frac{X-\mu}{t}\right) \tag{18}$$

with $h_\varepsilon$ being the inverse of the derivative $f_\varepsilon'$. An analogous formula applies to the upper-bound. For the $\alpha$ divergence, we have $h_\varepsilon\left(y\right)=\left(1+\frac{\alpha-1}{\alpha}\varepsilon y\right)_+^{\frac{1}{\alpha-1}}$ under the assumption of **Proposition 3**. The Radon–Nikodým derivative $Z$ of (18) vanishes for large $X$ if $\alpha>1$, while it is positive if $\alpha\in\left(0,1\right)$. Note that the optimal pair $\left(\mu,t\right)$ as the Lagrangian multipliers naturally gives the normalization $\mathbb{E}\left[Z\right]=1$ (e.g., see, **Proof of Proposition 2**).

### 2.6 The gamma case

We analysed the gamma case, the representative case of probability distributions having exactly one mode, where the probability measure $\mathbb{P}$ associates the probability density function (PDF) of the form

$$p(u) = \frac{1}{\Gamma(a)b^a} u^{a-1} \exp\left(-\frac{u}{b}\right), \ u > 0 \tag{19}$$

with parameters $a, b > 0$ and the gamma function $\Gamma$. A variety of environmental variables follow a gamma distribution, such as short-term rainfall intensity (Martinez-Villalobos and Neelin, 2019) and annual maximum flow discharge (Zhang et al., 2019). The gamma distribution also plays a crucial role in modelling long-memory stochastic differential equations of the superposition type. Here, a superposed stochastic differential equation is formally expressed by a superposition (i.e., integration with respect to a Lévy basis) of mutually-independent mean-reverting linear or affine stochastic differential equations having different reversion speeds (Barndorff-Nielsen and Stelzer, 2013; Leonte et al., 2023; Yoshioka et al., 2023a). The gamma distribution then appears as the PDF of the reversion speed, and the reciprocal moment parameterizes the moments and tails of the autocorrelation of the superposed process.

We analyse the existence of Orlicz regret and divergence risk measures against the moment $X = X(u) = u^{\gamma}$ ($\gamma \in \mathbb{R}$, $\gamma \neq 0$), where $u$ follows a gamma distribution. The primary interest will be on average ($\gamma = 1$), while this setting covers not only the classical moments $\gamma \in \mathbb{N}$ but also the reciprocal moment $\gamma = -1$ that may arise in some other applications like superposed processes, as discussed above. For $\alpha = 1$, the following **Proposition 6** holds true. Its proof is omitted because it uses a straightforward calculation considering the integrability of exponential and non-exponential parts separately.

***Proposition 6***

*Assume $\alpha = 1$. Then, it follows that*

$$\mathbb{E}\left[g_{\varepsilon}\left(\frac{X(u)-\mu}{t}\right)\right] = \frac{1}{\varepsilon}\left\{\int_0^{+\infty} \frac{1}{\Gamma(a)b^a} u^{a-1} \exp\left(-\frac{u}{b} + \varepsilon\frac{u^{\gamma}-\mu}{t}\right)\mathrm{d}u - 1\right\}, \tag{20}$$

*and there is a pair $(\mu, t) \in \mathbb{R} \times (0, +\infty)$ such that the right-hand side of (20) exists only if $0 < \gamma \leq 1$. Similarly, it follows that*

$$\mathbb{E}\left[g_{\varepsilon}\left(-\frac{X(u)-\mu}{t}\right)\right] = \frac{1}{\varepsilon}\left\{\int_0^{+\infty} \frac{1}{\Gamma(a)b^a} u^{a-1} \exp\left(-\frac{u}{b} - \varepsilon\frac{u^{\gamma}-\mu}{t}\right)\mathrm{d}u - 1\right\} \tag{21}$$

*and there is a pair $(\mu, t) \in \mathbb{R} \times (0, +\infty)$ such that the right-hand side of (21) exists for any $\gamma \in \mathbb{R}$.*

The proposition below covers cases of $\alpha \neq 1$. The proof is again by a straightforward calculation.

***Proposition 7***

*Assume $\alpha > 1$. Then, it follows that*

$$\mathbb{E}\left[g_{\varepsilon}\left(\frac{X(u)-\mu}{t}\right)\right]=\frac{1}{\varepsilon}\left\{\int_{0}^{+\infty}\frac{1}{\Gamma(a)b^{a}}u^{a-1}\exp\left(-\frac{u}{b}\right)\left(1+\frac{\alpha-1}{\alpha}\varepsilon\frac{u^{\gamma}-\mu}{t}\right)_{+}^{\frac{\alpha}{\alpha-1}}\mathrm{d}u-1\right\}, \tag{22}$$

*and there is a pair* $(\mu,t)\in\mathbb{R}\times(0,+\infty)$ *such that the right-hand side of (22) exists if* $\gamma>-a\frac{\alpha-1}{\alpha}$, *or if* $\gamma=-a\frac{\alpha-1}{\alpha}$ *and* $\varepsilon>0$ *is sufficiently small. Furthermore, it also follows that*

$$\mathbb{E}\left[g_{\varepsilon}\left(-\frac{X(u)-\mu}{t}\right)\right]=\frac{1}{\varepsilon}\left\{\int_{0}^{+\infty}\frac{1}{\Gamma(a)b^{a}}u^{a-1}\exp\left(-\frac{u}{b}\right)\left(1-\frac{\alpha-1}{\alpha}\varepsilon\frac{u^{\gamma}-\mu}{t}\right)_{+}^{\frac{\alpha}{\alpha-1}}\mathrm{d}u-1\right\}, \tag{23}$$

*and there is a pair* $(\mu,t)\in\mathbb{R}\times(0,+\infty)$ *such that the right-hand side of* (23) *exists for any* $\gamma\in\mathbb{R}$ *.*

*Assume* $\alpha\in(0,1)$ *.Then, it follows that there is a pair* $(\mu,t)\in\mathbb{R}\times(0,+\infty)$ *such that the right-hand side of* (23) *exists for any* $\gamma\in\mathbb{R}$ *.*

According to **Propositions 6-7**, the use of the $\alpha$ divergence is essential for estimating the upper-bound if $\gamma<0$, while both the Kullback–Leibler and $\alpha$ divergences can be utilized for evaluating the lower-bound. The estimation of the average ( $\gamma=1$ ) can be performed without restrictions for all $\alpha>1$. In addition, the existence of the integrals in **Propositions 6-7** is analysed by focusing on either small or large $u>0$, that is, tail behaviour. The same strategy applies to related PDFs such as the *q*-Weibull one (Sánchez et al., 2023) and generalized Laplace one (Jiang et al., 2016) having the specific tails.

***Remark 5*** In general, the domains of Orlicz regrets and, hence, the divergence risk measures are different, even if they use common $f$ . For the gamma case, both the upper- and lower-bounds of the moment $u^{\gamma}$ ( $\gamma>0$ ) can be evaluated by the $\alpha$ divergence with $\alpha>1$. In this case, the problem setting is consistent.

### 2.7 Numerical algorithms

We present a gradient descent method to compute $\bar{R}_{f,\varepsilon}$ and $\underline{R}_{f,\varepsilon}$. Set

$$G(\mu,t)=\mu+t\left(\varepsilon+\mathbb{E}\left[g\left(\frac{X-\mu}{t}\right)\right]\right)\text{ for }\mu\in\mathbb{R}\text{ and }t>0\ (t\to\varepsilon t\text{ in (9) and (12)}), \tag{24}$$

with which we have

$$\bar{R}_{f,\varepsilon}(X)=\inf_{\substack{\mu\in\mathbb{R}\\t>0}}G(\mu,t)\text{ and }\underline{R}_{f,\varepsilon}(X)=\sup_{\substack{\mu\in\mathbb{R}\\t>0}}G(\mu,-t). \tag{25}$$

Assume that minimising $(\mu,t)$ of $\bar{R}_{f,\varepsilon}$ and maximising $(\mu,t)$ of $\underline{R}_{f,\varepsilon}$ exist. Subsequently, we expect that the minimiser and maximiser can be computed using the long-time limits of the solutions to ordinary differential equations with a pseudo-time $\tau\geq 0$:

$$\frac{\mathrm{d}}{\mathrm{d}\tau}(\mu,t)=\left(-\frac{\partial}{\partial\mu}G(\mu,t),-\frac{1}{\varepsilon}\frac{\partial}{\partial t}G(\mu,t)\right) \tag{26}$$

and

$$\frac{\mathrm{d}}{\mathrm{d}\tau}(\mu,t)=\left(\frac{\partial}{\partial\mu}G(\mu,t),\frac{1}{\varepsilon}\frac{\partial}{\partial t}G(\mu,t)\right) \tag{27}$$

under the suitable initial conditions. The scaling $\varepsilon^{-1}$ has been empirically identified to facilitate convergence but does not affect the steady states of (26)–(27). Substituting the long-term limit of each equation into $G(\mu,t)$ or $G(\mu,-t)$ yields the desired $\bar{R}_{f,\varepsilon}$ or $\underline{R}_{f,\varepsilon}$, respectively.

In the subsequent section, each expectation appearing on the right-hand sides of (26) and (27) is approximated by the quantile discretization (Yoshioka et al., 2023a) with 8,192 degrees of freedom (See, **Appendix B**), which was observed to be sufficiently fine. We implemented a gradient descent to solve these differential equations using a fully explicit Euler method with a small pseudo-time increment of 0.5. The gradient descent is terminated if the absolute difference between the current and previous $\mu$ and $t$ becomes smaller than $10^{-10}$. The initial conditions are set as $(\mu,t)=(0,10)$.

***Remark 6*** We can also compute $V_{g_\varepsilon}$ and $W_{g_\varepsilon}$ in the same manner by considering the gradient descent with respect to $t$ by fixing $\mu=0$.

## 3. Application

The Orlicz regrets and associated divergence risk measures were applied to real environmental data.

### 3.1 Study site

The study site was in the downstream reach of the Hii River, a first-class river in the San-in area of Japan (**Figure 2**). This river has the brackish lakes as Ramsar sites, which are Lakes Shinji and Nakaumi, serving as crucial habitats for various aquatic species such as the Japanese basket clam *Corbicula japonica* (Derot et al., 2022) and migratory fishes (Takahara et al., 2023). The water quality of the Hii River, particularly its downstream reaches, including Otsu, has been measured by several researchers because Lake Shinji suffers from eutrophication caused by excessive inputs of phosphorus and nitrogen owing to intensive agricultural activities and runoff from forests in the watershed (See, **Figure 2**) (Somura et al., 2012; Ide et al., 2019). The threshold values of several key water quality indices, such as total nitrogen (TN) and total phosphorous (TP), for rivers and lakes are publicly provided, as discussed later in this paper. To the best of our knowledge, a stochastic analysis of water quality at the study site, including a water quality assessment under the model uncertainty on which we are focusing, has not been performed for the study site.

We obtained the unique 30-year weekly data for two indices from August 20, 1991, to December 27, 2021: TN and TP at the Otsu station from Dr. Ikuo Takeda, the author of (Takeda et al., 2023) (**Figure 3**). Both TN and TP are indicators of aquatic environments, such that high concentrations cause eutrophication (e.g., Jia et al., 2023; Phlips et al., 2023). By contrast, low concentrations of these water quality indices potentially result in the poor growth of riparian vegetation (Qian et al., 2023).

Otsu is located 12.3 (km) upstream of Lake Shinji, Japan. There is no tributary river between Otsu and Lake Shinji, although there are agricultural, forested, and residential areas surrounding and upstream of Otsu. Therefore, modeling and analysis of water quality dynamics at Otsu are important for assessing the water environment of the downstream part of the Hii River, particularly Lake Shinji.

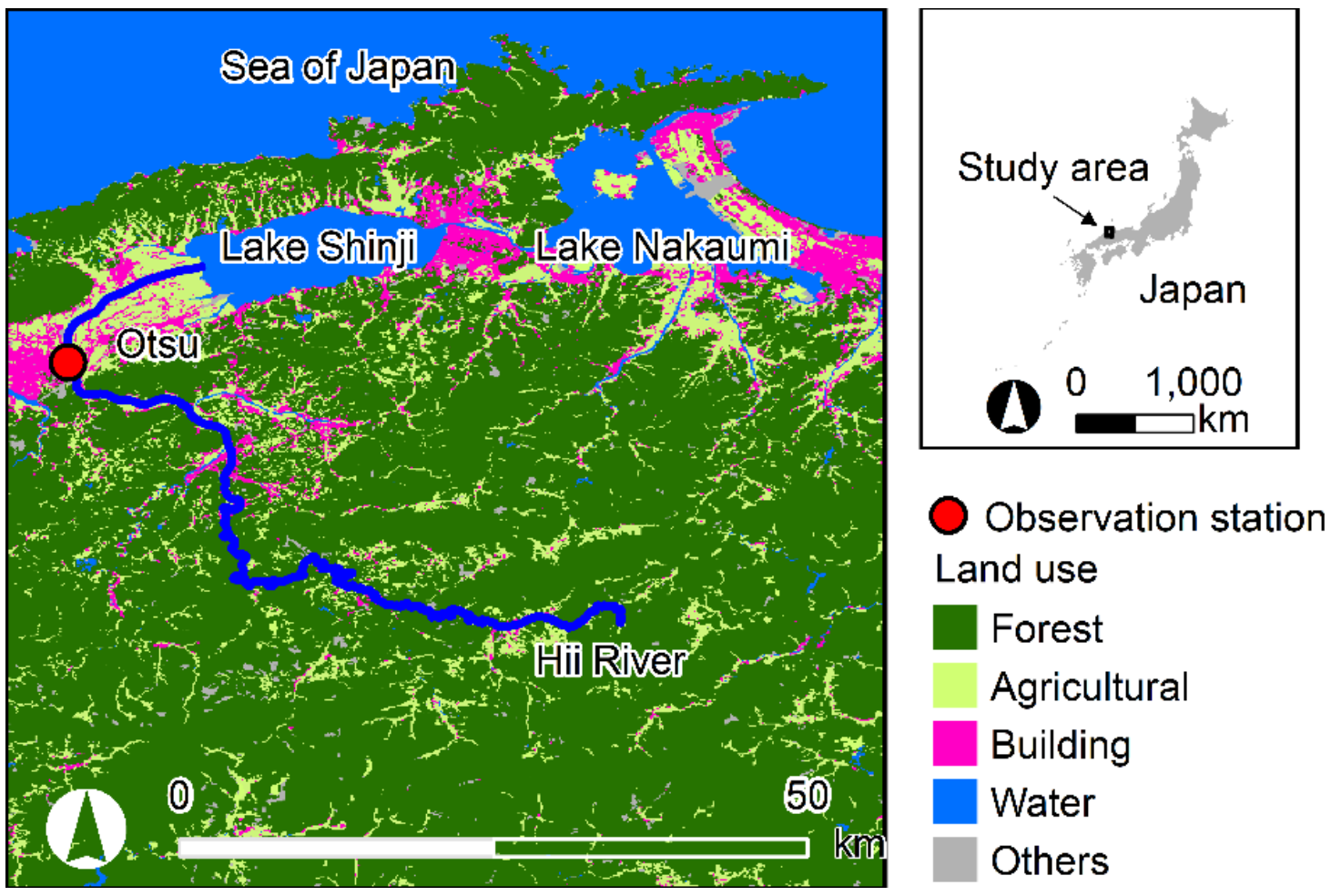

**Figure 2.** Map around the study site.

### 3.2 Model identification

As Takeda (2023) argued that water quality had not shown significant trends since 1991, we used all the data from August 1991 to December 2021 by assuming its stationarity. Let $\{t_k\}_{k=1,2,3,\ldots K}$ be a strictly increasing sequence representing the sampling date with the total amount of data $K \in \mathbb{N}$ and $K \geq 2$. The observations on day $t_k$ are denoted as $Y_{t_k}$. The observation value between the semi-closed time interval $[t_k, t_{k+1})$ was set as $Y_{t_k}$ ($k = 1, 2, 3, \ldots, K-1$). The units of the TN and TP, and the measurement methods used were as follows: TN (mg/L, ultraviolet absorption spectrophotometry after decomposition by potassium persulfate) and TP (mg/L, molybdenum blue method after decomposition by potassium persulfate).

We fit the gamma distribution (19) to TN and TP. Parameters $a, b$ were identified through the two equations obtained from the moment-matching method and are expressed as follows:

$$\text{Mean: } ab = \sum_{k=1}^{K} Y_{t_k} \text{ , Variance: } ab^2 = \left( \frac{1}{K-1} \sum_{k=1}^{K} \left( Y_{t_k} - \frac{1}{K} \sum_{l=1}^{K} Y_{t_l} \right)^2 \right). \quad (28)$$

**Table 1** lists the fitted parameter values for each water quality index obtained by applying the identification methods. **Table 1** shows that only TP satisfies $a < 1$, whereas TN satisfies $a > 1$. Therefore, the two indices had qualitatively different histogram shapes according to (19) and **Figure 4**. **Table 2** lists the empirical and theoretical moments of the fitted model of each index, showing that the average and variance were correctly captured by the proposed identification method, while the higher-order moments were less accurately reproduced. This is because the extremely high values present in the empirical data are not captured by the fitted model. This type of error motivates the approach to account for uncertainty. However, the signs of skewness and kurtosis were correctly captured by the fitted model, which is also shown in **Figure 4** when comparing the empirical and fitted histograms.

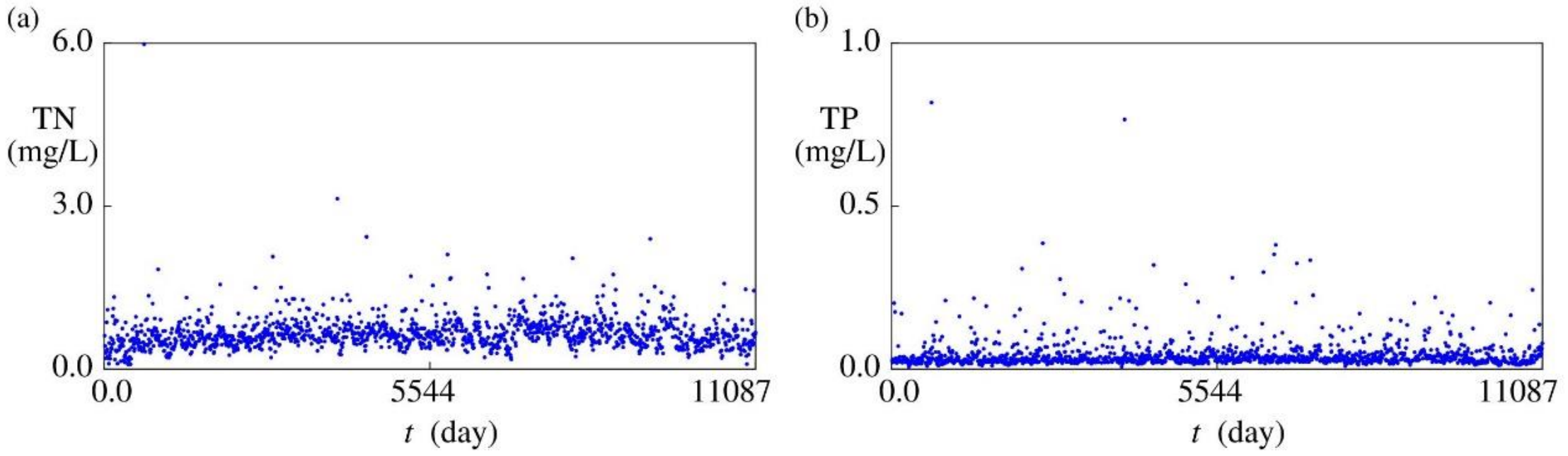

**Figure 3.** Time series data of (a) TN and (b) TP.

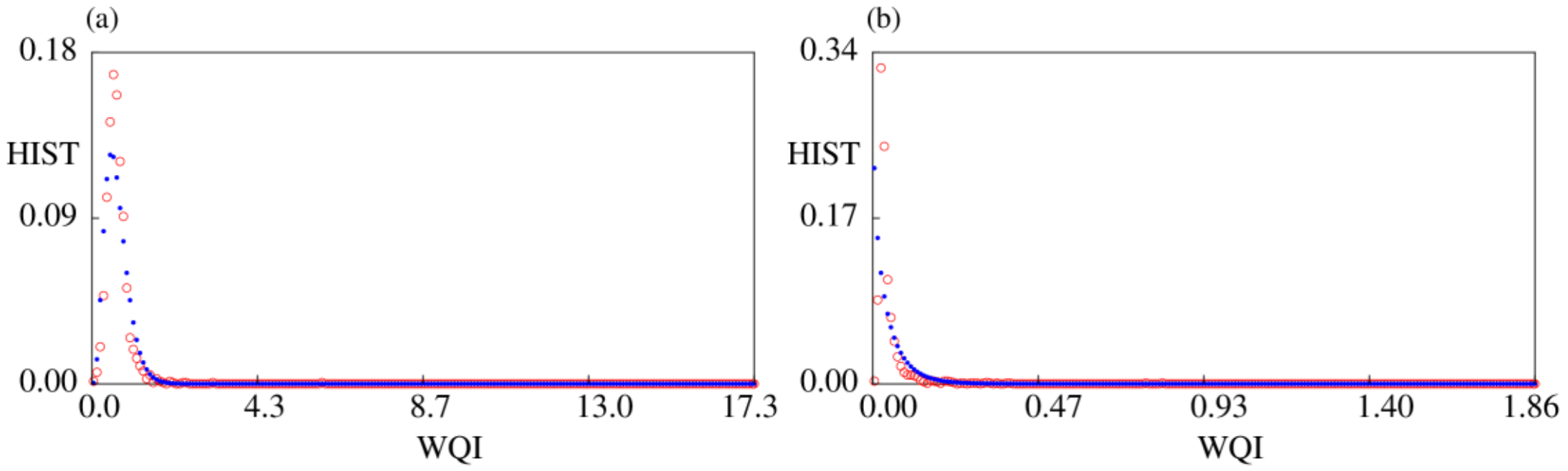

**Figure 4**. Comparison of empirical and fitted histograms (HIST) of (a) TN and (b) TP.

**Table 1.** Fitted parameter values.

| | TN | TP |
|---|---|---|
| $a$ (-) | 4.60.E+00 | 8.01.E-01 |
| $b$ (mg/L) | 1.42.E-01 | 5.41.E-02 |

**Table 2.** Comparison of empirical (Emp) moments and theoretical ones by the fitted model (Fit).

| | TN | | TP | |
|---|---|---|---|---|
| | Emp | Fit | Emp | Fit |
| Average (mg/L) | 6.53.E-01 | 6.53.E-01 | 4.34.E-02 | 4.34.E-02 |
| Variance ($mg^2/L^2$) | 9.26.E-02 | 9.26.E-02 | 2.35.E-03 | 2.35.E-03 |
| Skewness (-) | 5.10.E+00 | 9.33.E-01 | 7.39.E+00 | 2.23.E+00 |
| Kurtosis (-) | 7.05.E+01 | 1.31.E+00 | 8.85.E+01 | 7.49.E+00 |

### 3.3 Computational results

We computed the Orlicz regrets, divergence risk measures, and the corresponding Radon–Nikodým derivatives for TN and TP concentrations. The computational results below show how the target expectation of the water quality indices should be evaluated depending on the uncertainty aversion $\varepsilon$.

**Figures 5(a)-(b)** show the computed Orlicz regrets and divergence risk measures for TN and TP as functions of the uncertainty aversion parameter $\varepsilon = 10^{-5+6i/1000}$ ($i = 0,1,2,...,800$) with $\alpha = 1$. **Figures**

**6–7** show the computed Orlicz regrets and divergence risk measures for $\alpha = 0.5$ and $\alpha = 1.5$, respectively. The theoretical ordering of **Proposition 4** is satisfied for all the computations presented here. For both the upper- and lower-bounds, the case $\alpha = 1.5$ less strongly overestimates and underestimates the expectations than the case $\alpha = 1.0$; By contrast, the case $\alpha = 0.5$ more strongly underestimates the expectations. The α divergence with the parameter $\alpha$ larger than one (smaller than one) therefore leads to more optimistic (more pessimistic) estimates than the case $\alpha = 1.0$.

**Figures 8(a)-(b)** show the Radon–Nikodým derivatives to provide the worst-case estimations for TN and TP, where we assume $\alpha = 1.5$. For $\alpha = 0.5$, the gradient descent failed to converge against the upper bounds, which is considered the numerical appearance of the first statement in **Proposition 3**. For both TN and TP, increasing the uncertainty aversion $\varepsilon$ results in the PDFs $p$ being more weighted on a larger concentration ( concentration) for the upper-bounding case ( lower-bounding case). In particular, for the TN, the horizontal shift of the maximum point of the PDF according to uncertainty aversion is visible. By contrast, the unimodality of the PDF $p$ of the TP is preserved for the examined cases. These observations are qualitatively the same for different values of $\alpha$, except for the upper-bounding case with $\alpha \in (0,1)$, as implied above for $\alpha = 0.5$, where the expectation becomes ill-posed.

We finally analyse the threshold probability of the environmental criteria defined by the Ministry of Environment, Japan. For both TN and TP, there exists some prescribed threshold value $\bar{X} > 0$ above which the river water is not suited for the environmental conservation, agricultural use, and industrial use drinking purpose (https://www.env.go.jp/kijun/wt2-1-2.html). The threshold values $\bar{X}$ are 1.0 (mg/L) for TN and 0.1 (mg/L) for TP. Subsequently, we numerically compute the safety probability $P = \Pr\left(X \le \bar{X}\right) = \mathbb{E}_{\mathbb{Q}(\phi)}\left[\mathbb{I}\left(X \le \bar{X}\right)\right]$ under the worst-case uncertainty using both the lower-bounding (optimistic) and upper-bounding (pessimistic) cases in the context of the divergence risk measures. Higher values of the safety probability $P$ are preferred to avoid the eutrophication as an extreme event. Its value is 0.87 and 0.89 for TN and TP, respectively, when there is no uncertainty.

**Figure 9** shows the computed $P$ for different values of $\alpha$. In the absence of uncertainties, the safety probabilities are 0.13 and 0.11 for TN and TP, respectively. In these figures, the results corresponding to the upper-bounded cases (resp., lower-bounding cases) of the expectation $\mathbb{E}[X]$ appear below (above) the grey line standing for the safety probability without uncertainty. The upper-bounding case therefore leads to more pessimistic and hence, smaller safety probability, and vice versa. For both TN and TP, the cases with $\alpha = 1.5$ deviate less from the case $\alpha = 1.0$ without uncertainties being consistent with the results in **Figures 5–6**. Particularly, the case $\alpha = 1.0$ provides at maximum 10% and 7% smaller safety probability for TN and TP under the computational conditions, respectively. Finally, the overestimation of the safety probability is comparable for the different values of $\alpha$. This result is not so useful for TN and TP, while it can be useful for analysing the other water quality indices such as the dissolved oxygen that should not be too low (Wang et al., 2022) if the data will become available.

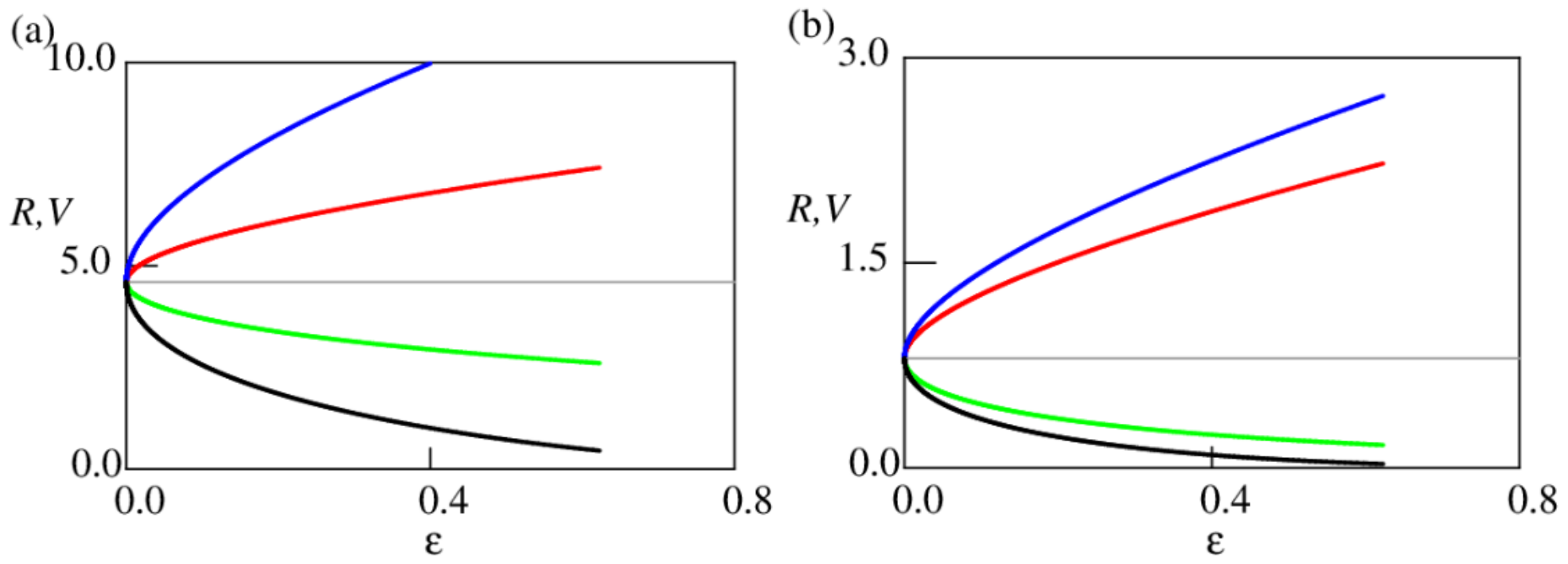

**Figure 5.** The Orlicz regrets and divergence risk measures as functions of $\varepsilon$: (a) TN and (b) TP ($\alpha = 1$). The colours represent $R = \overline{R}_{f,\varepsilon}$ (Red), $V = V_{g_\varepsilon}$ (Blue), $R = \underline{R}_{f,\varepsilon}$ (Green), and $V = W_{g_\varepsilon}$ (Black).

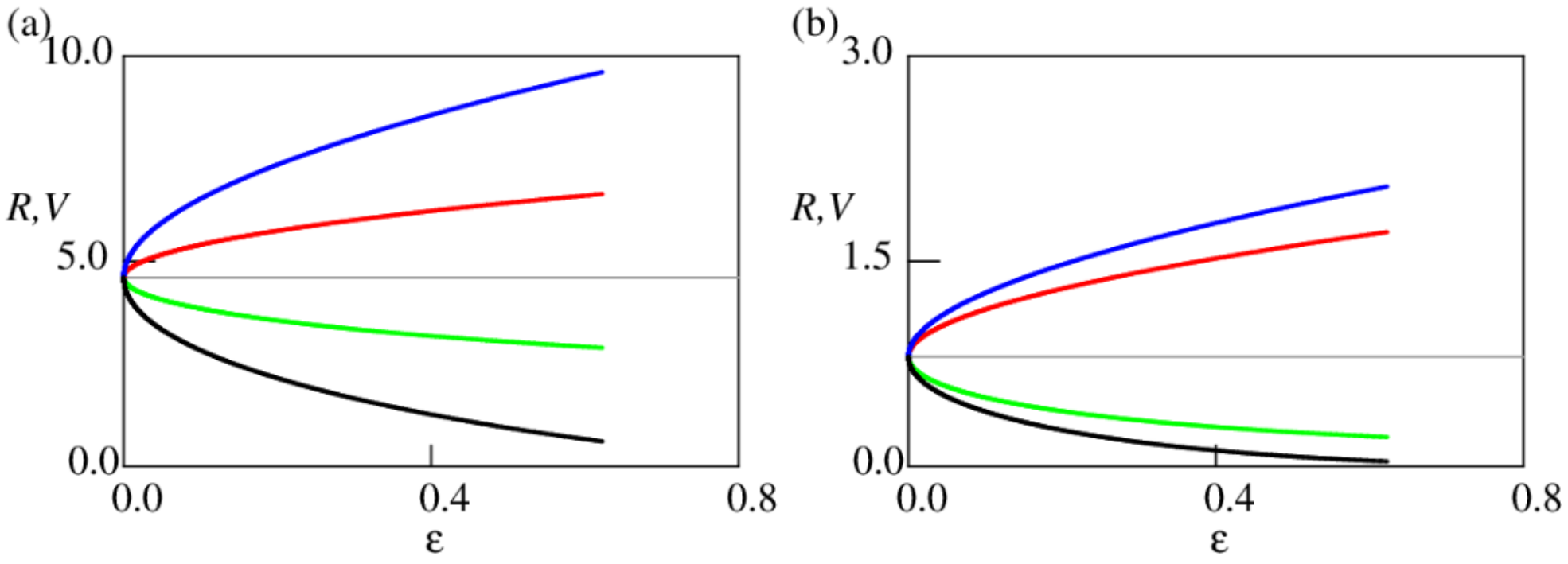

**Figure 6.** The Orlicz regrets and divergence risk measures for $\alpha = 1.5$: (a) TN and (b) TP. The legends are the same with **Figure 5**.

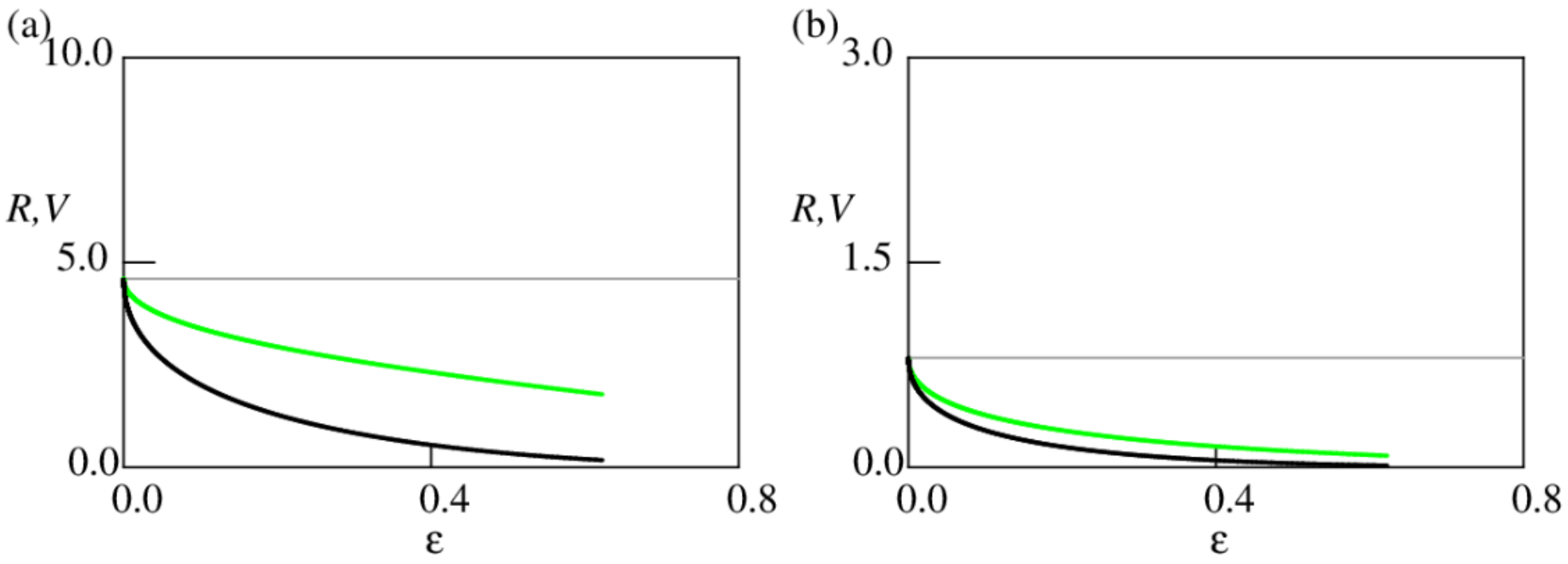

**Figure 7.** The Orlicz regrets and divergence risk measures for $\alpha = 0.5$: (a) TN and (b) TP. The legends are the same with **Figure 5**.

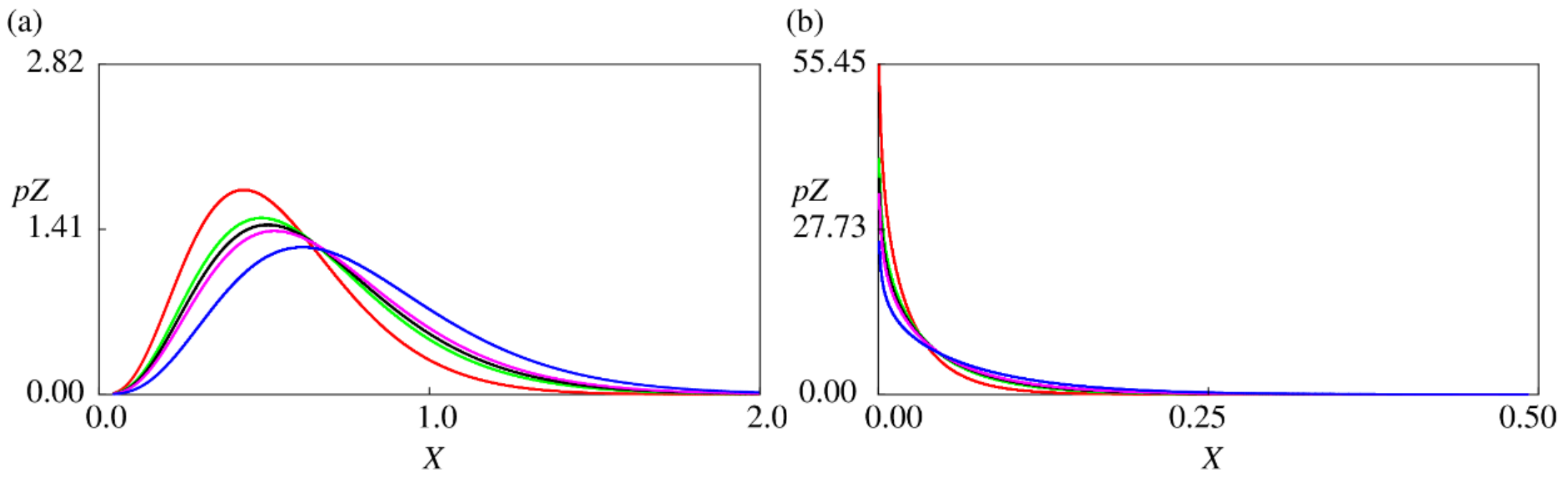

**Figure 8.** The distorted PDF $pZ$ for (a) TN and (b) TP ($\alpha = 1.5$). The colours represent the lower-bounding case with $\varepsilon = 0.1$ (Red), lower-bounding case with $\varepsilon = 0.005$ (Green), no uncertainty (Black), upper-bounding case with $\varepsilon = 0.005$ (Pink), and upper-bounding case with $\varepsilon = 0.1$ (Blue).

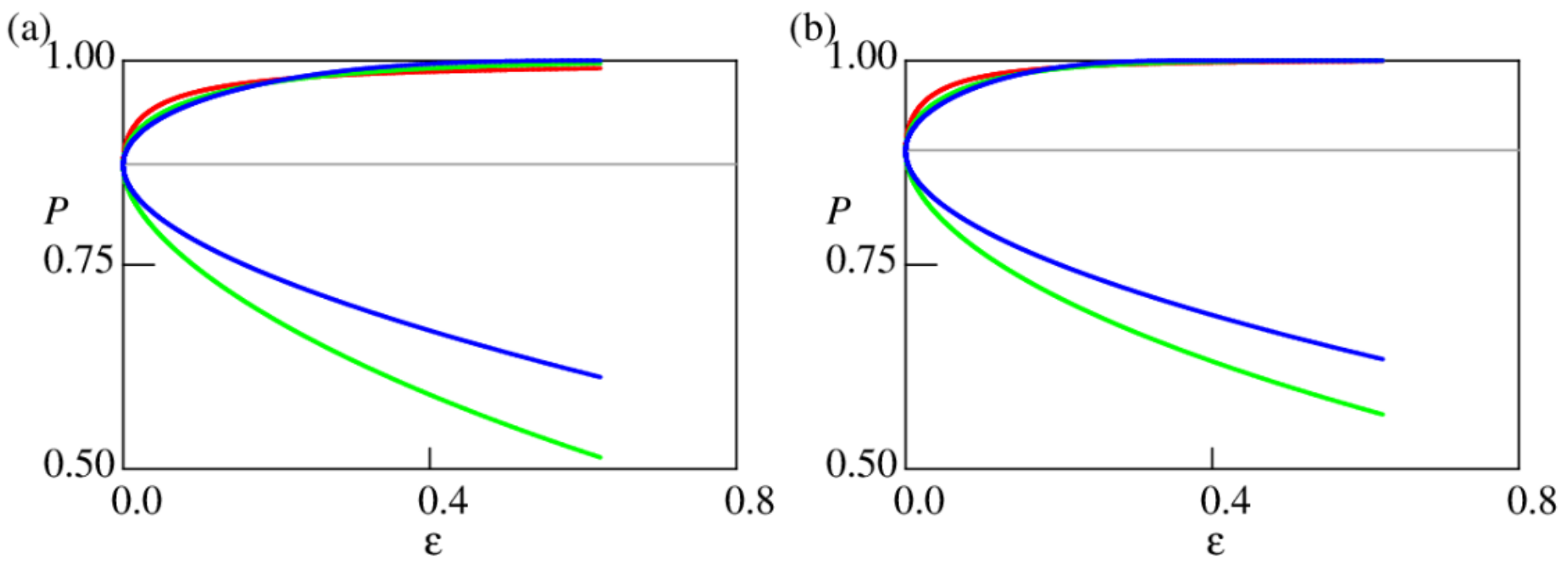

**Figure 9.** The safety probability $P$ under the worst-case estimations for (a) TN and (b) TP. The grey line stands for the safety probability under no uncertainty. The plots above the grey line correspond to $P$ for the lower-bounds of $\mathbb{E}[X]$ that are more optimistic, while the plots below the grey line correspond to $P$ for the upper-bounds of $\mathbb{E}[X]$ that are more pessimistic. The colours represent $\alpha = 0.5$ (Red), $\alpha = 1.0$ (Green), and $\alpha = 1.5$ (Blue).

## 4. Conclusion

A pair of Orlicz regrets was proposed and their properties were mathematically analysed. Numerical algorithms for computing them were also proposed. Finally, they were applied to the real water environmental data from Japanese river environments using the gradient descent-type method. In this study, the proposed Orlicz regrets and the divergence risk measures were consistently established within a unified mathematical framework. Particularly, in conventional studies, Young's and entropic penalties are defined separately, while our framework can simultaneously specify them based on the divergence.

The Orlicz regrets can be applied not only to the environmental variables, but also to the other random variables arising in a wide variety of problems; potential examples include but are not limited to the portfolio optimization (Das et al., 2017), commodity prices (Li et al., 2023), and guarantees of system performance in general (Salamati and Zamani, 2022). This study focuses on static problems, although dynamic problems can also be defined using recursive equations to sequentially optimize the target goal (Tian et al., 2023; Gu et al., 2023). Exploring the relationship between the maximum entropy principle and Orlicz regrets and divergence risk measures will also be interesting (Villa-Morales and Rincón, 2023). The authors will investigate Orlicz regrets that are suitable for streamflow discharge in mountainous river environments as jump-driven stochastic processes.

## Appendices of “Orlicz regrets to consistently bound statistics of random variables with an application to environmental indicators”

### Appendix A: Proofs of Propositions

***Proof of Proposition 1***

The second statement is firstly proven We have the equality

$$\begin{aligned} W_g(X) &= \sup_{t>0}\left\{-t\left(1+\mathbb{E}\left[g\left(\frac{-X}{t}\right)\right]\right)\right\} \\ &= \sup_{t>0}\left\{-t\left(1+\sup_{Z\geq 0,\,\mathbb{E}[Z]<+\infty}\left(\mathbb{E}\left[\frac{-X}{t}Z\right]-\mathbb{E}\left[f(Z)\right]\right)\right)\right\} \\ &= \sup_{t>0}\left\{-t+t\inf_{Z\geq 0,\,\mathbb{E}[Z]<+\infty}\left(\mathbb{E}\left[\frac{X}{t}Z\right]+\mathbb{E}\left[f(Z)\right]\right)\right\} \\ &= \sup_{t>0}\inf_{Z\geq 0,\,\mathbb{E}[Z]<+\infty}\left\{\mathbb{E}[XZ]+t\left(\mathbb{E}\left[f(Z)\right]-1\right)\right\} \end{aligned}. \tag{29}$$

The remaining part follows by the classical Lagrangian representation of a concave minimisation problem

Now, we prove **W1-W5** in this order

**W1** We have, by the transformation use $t \to \lambda s$,

$$W_g(\lambda X) = \sup_{t>0}\left\{-t\left(1+\mathbb{E}\left[g\left(\frac{-\lambda X}{t}\right)\right]\right)\right\} = \lambda\sup_{s>0}\left\{-s\left(1+\mathbb{E}\left[g\left(\frac{-X}{s}\right)\right]\right)\right\} = \lambda W_g(X). \tag{30}$$

**W2** For any $t_1, t_2 > 0$, we have the equality

$$W_g(X+Y) = \sup_{t>0}\left\{-t\left(1+\mathbb{E}\left[g\left(\frac{-(X+Y)}{t}\right)\right]\right)\right\} \geq -\frac{t_1+t_2}{t_1t_2}\left(1+\mathbb{E}\left[g\left(-(X+Y)\frac{t_1t_2}{t_1+t_2}\right)\right]\right). \tag{31}$$

We then obtain

$$\begin{aligned} W_g(X+Y) &\geq -\frac{t_1+t_2}{t_1t_2}\left(1+\mathbb{E}\left[g\left(-\frac{t_2}{t_1+t_2}t_1X-\frac{t_1}{t_1+t_2}t_2Y\right)\right]\right) \\ &\geq -\frac{t_1+t_2}{t_1t_2}\left(1+\frac{t_2}{t_1+t_2}\mathbb{E}\left[g(-t_1X)\right]+\frac{t_1}{t_1+t_2}\mathbb{E}\left[g(-t_2Y)\right]\right), \\ &= -\frac{1}{t_1}\left(1+\mathbb{E}\left[g(-t_1X)\right]\right)-\frac{1}{t_2}\left(1+\mathbb{E}\left[g(-t_2X)\right]\right) \end{aligned} \tag{32}$$

where we used the convexity of $g$ to obtain the second line of (32) Because $t_1, t_2 > 0$ are arbitrary, considering the supremum with respect to them in (32) proves **W2.**

**W3** This is owing to the increasing property of $g$ For any $t > 0$, we obtain

$$\mathbb{E}\left[g\left(\frac{-X}{t}\right)\right] \geq \mathbb{E}\left[g\left(\frac{-Y}{t}\right)\right] \text{ and thus } -t\left(1+\mathbb{E}\left[g\left(\frac{-\lambda X}{t}\right)\right]\right) \leq -t\left(1+\mathbb{E}\left[g\left(\frac{-\lambda Y}{t}\right)\right]\right). \tag{33}$$

The desired monotonicity is proven by taking the supremum with respect to $t > 0$ in (32)

**W4** For any $y \in \mathbb{R}$, we have

$$g(y) = \sup_{x \in \mathbb{R}} \{xy - f(x)\} \geq y - f(1) = y \,. \tag{34}$$

We also have

$$\mathbb{E}[X] = -t\mathbb{E}\left[\frac{-X}{t}\right] \geq -t - t\mathbb{E}\left[\frac{-X}{t}\right] = -t\left(1 + \mathbb{E}\left[\frac{-X}{t}\right]\right),\ t > 0 \,. \tag{35}$$

Considering (34), (35) leads to

$$\mathbb{E}[X] \geq -t\left(1 + \mathbb{E}\left[\frac{-X}{t}\right]\right) \geq -t\left(1 + \mathbb{E}\left[g\left(\frac{-X}{t}\right)\right]\right),\ t > 0 \,. \tag{36}$$

Taking the supremum with respect to $t > 0$ in (36) proves **W4**

**W5** Assume that almost surely $X = 0$ Then, because of $g(0) = 0$, it follows that

$$W_g(0) = \sup_{t>0}\left\{-t\left(1 + \mathbb{E}\left[g\left(\frac{0}{t}\right)\right]\right)\right\} = \sup_{t>0}\left\{-t\left(1 + \mathbb{E}\left[g(0)\right]\right)\right\} = 0 \,. \tag{37}$$

Conversely, assume that the $W_g(X) = 0$ for some $X \leq 0$ not almost surely $X = 0$ We have $\mathbb{E}[X] < 0$ Then, by **V4**, we obtain

$$0 > \mathbb{E}[X] \geq W_g(X) \,. \tag{38}$$

Consequently, we obtain the desired point-separating property $0 > W_g(X)$

□

***Proof of Proposition 2***

First, $\underline{R}_{f,\varepsilon}$ is rewritten as follows:

$$\begin{aligned}\underline{R}_{f,\varepsilon}(X) &= \inf\left\{\mathbb{E}_{\mathbb{P}}\left[\frac{\mathrm{d}\mathbb{Q}}{\mathrm{d}\mathbb{P}}X\right] \,\middle|\, \mathbb{Q} \in \mathfrak{Q},\ \mathbb{E}_{\mathbb{P}}\left[f_\varepsilon\left(\frac{\mathrm{d}\mathbb{Q}}{\mathrm{d}\mathbb{P}}\right)\right] \leq 1\right\} \\ &= \inf\left\{\mathbb{E}_{\mathbb{P}}[ZX] \,\middle|\, Z \geq 0, \mathbb{E}_{\mathbb{P}}[Z] = 1,\ \mathbb{E}_{\mathbb{P}}\left[f_\varepsilon(Z)\right] \leq 1\right\} \\ &= -\sup\left\{\mathbb{E}_{\mathbb{P}}[-XZ] \,\middle|\, Z \geq 0, \mathbb{E}_{\mathbb{P}}[Z] = 1,\ \mathbb{E}_{\mathbb{P}}\left[f_\varepsilon(Z)\right] \leq 1\right\}\end{aligned} \,. \tag{39}$$

The classical method of Lagrangian multipliers leads to

$$\underline{R}_{f,\varepsilon}(X) = -\inf_{\substack{\mu \in \mathbb{R} \\ t>0}} \sup_{Z \geq 0}\left\{\mathbb{E}_{\mathbb{P}}[-XZ] - \mu\left(\mathbb{E}_{\mathbb{P}}[Z] - 1\right) - t\left(\mathbb{E}_{\mathbb{P}}\left[f_\varepsilon(Z)\right] - 1\right)\right\}, \tag{40}$$

where the supremum is taken with respect to non-negative random variables $Z$ We continue as follows:

$$
\begin{aligned}
\underline{R}_{f,\varepsilon}(X) &= -\inf_{\substack{\mu\in\mathbb{R}\\ t>0}} \sup_{Z\ge 0}\left\{\mu + t + \mathbb{E}_{\mathbb{P}}[-XZ] - \mu\mathbb{E}_{\mathbb{P}}[Z] - t\mathbb{E}_{\mathbb{P}}\left[f_\varepsilon(Z)\right]\right\} \\
&= -\inf_{\substack{\mu\in\mathbb{R}\\ t>0}} \sup_{Z\ge 0}\left\{\mu + t + \mathbb{E}_{\mathbb{P}}\left[-(X+\mu)Z\right] - t\mathbb{E}_{\mathbb{P}}\left[f_\varepsilon(Z)\right]\right\} \\
&= -\inf_{\substack{\mu\in\mathbb{R}\\ t>0}} \sup_{Z\ge 0}\left\{\mu + t + t\left\{\mathbb{E}_{\mathbb{P}}\left[-\frac{X+\mu}{t}Z\right] - \mathbb{E}_{\mathbb{P}}\left[f_\varepsilon(Z)\right]\right\}\right\} \\
&= -\inf_{\substack{\mu\in\mathbb{R}\\ t>0}} \left\{\mu + t + t\sup_{Z\ge 0}\left\{\mathbb{E}_{\mathbb{P}}\left[-\frac{X+\mu}{t}Z\right] - \mathbb{E}_{\mathbb{P}}\left[f_\varepsilon(Z)\right]\right\}\right\} \\
&= -\inf_{\substack{\mu\in\mathbb{R}\\ t>0}} \left\{\mu + t + t\mathbb{E}_{\mathbb{P}}\left[g_\varepsilon\left(-\frac{X+\mu}{t}\right)\right]\right\} \\
&= -\inf_{\substack{\mu\in\mathbb{R}\\ t>0}} \left\{\mu + t\left(1 + \mathbb{E}_{\mathbb{P}}\left[g_\varepsilon\left(-\frac{X+\mu}{t}\right)\right]\right)\right\}
\end{aligned} \tag{41}
$$

We further continue to obtain the desired result as follows:

$$
\begin{aligned}
\underline{R}_{f,\varepsilon}(X) &= -\inf_{\substack{\mu\in\mathbb{R}\\ t>0}} \left\{(-\mu) + t\left(1 + \mathbb{E}_{\mathbb{P}}\left[g_\varepsilon\left(-\frac{X+(-\mu)}{t}\right)\right]\right)\right\} \\
&= -\inf_{\substack{\mu\in\mathbb{R}\\ t>0}} \left\{-\mu + t\left(1 + \mathbb{E}_{\mathbb{P}}\left[g_\varepsilon\left(-\frac{X-\mu}{t}\right)\right]\right)\right\} \\
&= \sup_{\substack{\mu\in\mathbb{R}\\ t>0}} \left\{\mu - t\left(1 + \mathbb{E}_{\mathbb{P}}\left[g_\varepsilon\left(-\frac{X-\mu}{t}\right)\right]\right)\right\} \\
&= \sup_{\mu\in\mathbb{R}} \left\{\mu + \sup_{t>0}\left[-t\left(1 + \mathbb{E}_{\mathbb{P}}\left[g_\varepsilon\left(-\frac{X-\mu}{t}\right)\right]\right)\right]\right\} \\
&= \sup_{\mu\in\mathbb{R}} \left\{\mu + W_{g_\varepsilon}(X-\mu)\right\}
\end{aligned} \tag{42}
$$

□

***Proof of Proposition 3***

First, we prove the statement of $\bar{R}_{f,\varepsilon}$ Assume that $\alpha\in(0,1)$ We must ensure that any minimising pair $(\mu,t)$ of $\bar{R}_{f,\varepsilon}$, if it exists, satisfies

$$
1+\frac{\alpha-1}{\alpha}\varepsilon\left(\frac{X-\mu}{t}\right)>0 \text{ or equivalently } 1>\frac{1-\alpha}{\alpha}\varepsilon\frac{X-\mu}{t}. \tag{43}
$$

However, this is impossible owing to the unboundedness assumption of $X$

We prove the statement of $\underline{R}_{f,\varepsilon}$ We must show that a maximising pair $(\mu,t)$ of $\underline{R}_{f,\varepsilon}$, if it exists, satisfies

$$
1+\frac{\alpha-1}{\alpha}\varepsilon\left(-\frac{X-\mu}{t}\right)>0 \text{ or equivalently } 1+\frac{1-\alpha}{\alpha}\varepsilon\frac{X-\mu}{t}>0. \tag{44}
$$

This is satisfied when

$$
1+\frac{1-\alpha}{\alpha}\varepsilon\left(\frac{0-\mu}{t}\right)>0 \text{ or equivalently } \frac{\alpha}{(1-\alpha)\varepsilon}>\frac{\mu}{t}. \tag{45}
$$

Assume that (44) is not satisfied, then we have

$$G(\mu,-t)=\mu-t\left(1+\mathbb{E}\left[g_{\varepsilon}\left(-\frac{X-\mu}{t}\right)\right]\right)=-\infty \tag{46}$$

owing to the blow up of $g_{\varepsilon}$, while we also have

$$G(0,-t)=0-t\left(1+\mathbb{E}\left[g_{\varepsilon}\left(-\frac{X}{t}\right)\right]\right)>-\infty \text{ for any } t>0. \tag{47}$$

Hence, any pair $(\mu,t)\in\mathbb{R}\times(0,+\infty)$ violating (44) does not maximise $G(\mu,t)$

□

***Proof of Proposition 4***

We have

$$\overline{R}_{f,\varepsilon}(X)=\inf_{\mu\in\mathbb{R}}\left\{\mu+V_{g_{\varepsilon}}(X-\mu)\right\}\leq 0+V_{g_{\varepsilon}}(X-0)=V_{g_{\varepsilon}}(X) \tag{48}$$

as well as

$$\underline{R}_{f,\varepsilon}(X)=\sup_{\mu\in\mathbb{R}}\left\{\mu+W_{g_{\varepsilon}}(X-\mu)\right\}\geq 0+W_{g_{\varepsilon}}(X-0)=W_{g_{\varepsilon}}(X). \tag{49}$$

These results combined with **V4** and **W4** obtain the desired inequality

□

***Proof of Proposition 5***

The inequalities in (15) are clear according to Definitions (8) and (11) The right inequality in (16) follows from

$$\begin{aligned}
W_{g_{\varepsilon_2}}(X)&=\sup_{t>0}\left\{-t\left(1+\mathbb{E}\left[g_{\varepsilon_2}\left(\frac{-X}{t}\right)\right]\right)\right\}\\
&=\sup_{t>0}\left\{-t\left(1+\mathbb{E}\left[\frac{1}{\varepsilon_2}g\left(\frac{-\varepsilon_2 X}{t}\right)\right]\right)\right\}\\
&=\sup_{t>0}\left\{-\frac{t}{\varepsilon_2}\left(\varepsilon_2+\mathbb{E}\left[g\left(\frac{-\varepsilon_2 X}{t}\right)\right]\right)\right\}\\
&=\sup_{t>0}\left\{-t\left(\varepsilon_2+\mathbb{E}\left[g\left(\frac{-X}{t}\right)\right]\right)\right\} \qquad \left(\frac{t}{\varepsilon_2}\to t\right)\\
&\leq\sup_{t>0}\left\{-t\left(\varepsilon_1+\mathbb{E}\left[g\left(\frac{-X}{t}\right)\right]\right)\right\}\\
&=W_{g_{\varepsilon_1}}(X)
\end{aligned}. \tag{50}$$

The left inequality of (16) follows in the same manner The left inequality of (17) follows from Proposition 314 of Fröhlich and Williamson (2023) Finally, the right inequality of (17) follows from

$$
\begin{aligned}
W_{g_{\varepsilon_2}}(X) &= \sup_{t>0}\left\{-t\left(\varepsilon_2+\mathbb{E}\left[g\left(\frac{-X}{t}\right)\right]\right)\right\} \\
&\geq \sup_{t>0}\left\{-t\left(\varepsilon_2+\frac{\varepsilon_2}{\varepsilon_1}\mathbb{E}\left[g\left(\frac{-X}{t}\right)\right]\right)\right\} \\
&= \sup_{t>0}\left\{-\left(\frac{\varepsilon_2}{\varepsilon_1}t\right)\left(\varepsilon_1+\mathbb{E}\left[g\left(\frac{-X}{t}\right)\right]\right)\right\} \\
&= \frac{\varepsilon_2}{\varepsilon_1}W_{g_{\varepsilon_1}}(X)
\end{aligned}, \tag{51}
$$

where we used $\frac{\varepsilon_2}{\varepsilon_1} \geq 1$ and $-g_{\varepsilon_1}\left(\frac{-X}{t}\right) \leq 0$ for $X \leq 0$ The latter follows from the fact that $f$ is twice continuously differentiable and convex, and hence the first-order derivative $f'$ is increasing for $x > 0$, as well as $f(x) = +\infty$ for any $x < 0$ Indeed, for any $y \geq 0$, we have

$$
g(y) = \sup_{x\in\mathbb{R}}\{xy - f(x)\} = \sup_{x\geq 0}\{xy - f(x)\}. \tag{52}
$$

This is because for any $x < 0$, we have $xy - f(x) = -\infty$, whereas $0y - f(0) = -f(0) > -\infty$ From (52), for any $y_1, y_2 \in \mathbb{R}$ with $y_1 \leq y_2$, we have

$$
g(y_1) = \sup_{x\geq 0}\{xy_1 - f(x)\} \leq \sup_{x\geq 0}\{xy_2 - f(x)\} = g(y_2), \tag{53}
$$

from which we obtain

$$
g(y) \geq g(0) = \sup_{x\geq 0}\{-f(x)\} = 0 \text{ for any } y \geq 0, \tag{54}
$$

and hence $-g_{\varepsilon_1}\left(\frac{-X}{t}\right) \leq 0$ for $X \leq 0$

□

**Appendix B**

In order to be self-contained, we explain the quantile discretization for the gamma-type $p$ in Yoshioka et al (2023a) to approximate its distribution function Set $m \in \mathbb{N}$ and $N = 2^m$ For each $n \in \mathbb{N}$, there is the quantile $\alpha_{i,n}$ such that

$$
\frac{i}{n} = \int_0^{\alpha_{i,n}} p(x)\,\mathrm{d}x,\ i = 0,1,2,3,\ldots,n, \tag{55}
$$

where we formally set $\alpha_{i,n} = +\infty$ Note that $\int_{\alpha_{2i-2,2N}}^{\alpha_{2i,2N}} p(y)\,\mathrm{d}y = N^{-1}$ ($i = 1,2,3,\ldots,N$).

The Markovian lift with the degree-of-freedom $N\left(=2^m\right)$ approximates the distribution function $F(x) = \int_0^x p(y)\,\mathrm{d}y$ by $F_N(x) = \int_0^x \frac{1}{N}\sum_{i=1}^{N}\delta_{\{y=\alpha_{2i-1,2N}\}}\mathrm{d}y$, where $\delta_{\{y=\cdot\}}$ represents the Dirac's delta This is indeed a uniform approximation because, for each $x \in \left[\alpha_{2i-2,2N}, \alpha_{2i,2N}\right)$ ($i = 1,2,3,\ldots,N$), we have

$$
\begin{aligned}
\left|F(x)-F_N(x)\right| &= \left|\int_0^x p(y)\mathrm{d}y - \int_0^x N^{-1}\sum_{j=1}^{N}\delta_{\{y=\alpha_{2j-1,2N}\}}\mathrm{d}y\right| \\
&= \left|\int_{\alpha_{2i-2,2N}}^x p(y)\mathrm{d}y - \int_{\alpha_{2i-2,2N}}^x N^{-1}\sum_{i=1}^{N}\delta_{\{y=\alpha_{2j-1,2N}\}}\mathrm{d}y\right| \\
&= \left|\int_{\alpha_{2i-2,2N}}^x p(y)\mathrm{d}y - N^{-1}\int_{\alpha_{2i-2,2N}}^x \delta_{\{y=\alpha_{2i-1,2N}\}}\mathrm{d}y\right| \\
&\le \max\left\{\int_{\alpha_{2i-2,2N}}^x p(y)\mathrm{d}y, N^{-1}\int_{\alpha_{2i-2,2N}}^x \delta_{\{y=\alpha_{2i-1,2N}\}}\mathrm{d}y\right\} \\
&\le \max\left\{\int_{\alpha_{2i-2,2N}}^{\alpha_{2i,2N}} p(y)\mathrm{d}y, N^{-1}\int_{\alpha_{2i-2,2N}}^{\alpha_{2i,2N}} \delta_{\{y=\alpha_{2i-1,2N}\}}\mathrm{d}y\right\} \\
&\le N^{-1}
\end{aligned}
\quad , \qquad (56)
$$

where we used the continuity of $p$ and $\int_0^{\alpha_{2i-2,2N}} p(y)\mathrm{d}y = \int_0^{\alpha_{2i-2,2N}} \frac{1}{N}\sum_{j=1}^{N}\delta_{\{y=\alpha_{2j-1,2N}\}}\mathrm{d}y$ ( $i=1,2,3,...,N$ ) By Theorem 1316 and 1323 of Klenke (2020), $F_n$ weakly converges to $F$ The convergence speed is first-order in $N$ because of the estimate (56) The expectation of function $q(x)$ is then naively set as $N^{-1}\sum_{j=1}^{N} q(\alpha_{2j-1,2N})$ Due to the relationship $\int_{\alpha_{2i-2,2N}}^{\alpha_{2i,2N}} p(y)\mathrm{d}y = N^{-1}$ , the Markovian lift is seen as a deterministic importance sampling method for probability distributions In practice, each $\alpha_{i,n}$ can be found through some common numerical method such as a bisection method combined with a numerical or analytical quadrature formula